     \documentclass[11pt]{amsart}

\newcommand{\lrw}{\longrightarrow}

\DeclareMathOperator{\coz}{{coz}}
\DeclareMathOperator{\supp}{{supp}}

\DeclareMathOperator{\SOT}{SOT}

     \newtheorem{thm}{Theorem}

  \theoremstyle{definition}

	\newtheorem*{rem-eg}{Remark and Example}  

  \theoremstyle{remark}
\begin{document}

\parskip=0.3\baselineskip
\baselineskip=1.2\baselineskip

\title{Into isometries of $C_0(X,E)$'s}


\author{Jyh-Shyang Jeang}
\address{{Department of Applied Mathematics}\\
{National Sun Yat-sen University}\\
{Kaohsiung, 80424, Taiwan, R.O.C.}}
\email{jeangjs@math.nsysu.edu.tw}

\author{Ngai-Ching Wong}
\address{{Department of Applied Mathematics}\\
{National Sun Yat-sen University}\\
{Kaohsiung, 80424, Taiwan, R.O.C.}}
\email{wong@math.nsysu.edu.tw}

\thanks{Partially supported by 
National Science Council of Republic of China
(NSC 85-2121-M-110-006)}

\keywords{
Banach-Stone Theorem, into isometries, strictly convex Banach  
spaces, vector-valued 
functions}

\subjclass{46}

\begin{abstract}
Suppose  $X$ and $Y$ are locally  compact Hausdorff  spaces, $E$  
and $F$ are
Banach spaces and $F$ is strictly convex. We show that every linear  
isometry
$T$  from  $C_0(X,E)$  {\em  into}  $C_0(Y,F)$  is  essentially  a   
weighted
composition operator $Tf(y) = h(y) (f(\varphi(y)))$.
\end{abstract}

\maketitle

    Let $X$ and $Y$ be locally compact Hausdorff spaces, $E$ and  
$F$ Banach spaces and
$C_0(X,E)$ and $C_0(Y,F)$ the Banach  spaces  of continuous   
$E$-valued  and
$F$-valued   functions  defined  on  $X$  and  $Y$  vanishing  at   
infinity,
respectively.  Recall that a Banach space $E$ is said to be 
{\em strictly  convex}
if every  norm one element  of $E$ is an extreme  point  of the  
closed  unit
ball $U_E$ of $E$.  In \cite{Jer}, Jerison  gave a vector version of
Banach-Stone Theorem: If $X$ and $Y$ are compact Hausdorff spaces  
and $E$ is
a strictly  convex Banach
 space then every {\em surjective}  isometry  $T$ from $C(X,E)$
$(=C_0(X,E))$  onto  $C(Y,E)$  $(=C_0(Y,E))$  can  be written  as a 
{\em weighted
composition  operator}, i.e, $Tf(y)  = h(y) (f(\varphi(y)))$,  
$\forall  f \in
C(X,E)$, $\forall y \in Y$, where $\varphi$ is a homeomorphism from  
$Y$ onto
$X$ and $h$ is a continuous  map from  $Y$ into the space   
($B(E,E)$,$\SOT$) of
bounded linear operators from $E$ into $E$ equipped with the strong  
operator
topology
($\SOT$)  such that $h(y)$ is an isometrically  isomorphism  from  
$E$ onto $E$
for all $y$ in $Y$.  After then several generalizations  of Banach-Stone
Theorem in this direction have 
appeared  (see, for example, 
\cite{Beh}). We shall show in this note:

\medskip
\begin{thm}\label{thm}
Suppose  $X$ and $Y$ are locally  compact Hausdorff  spaces, $E$  
and $F$ are
Banach spaces, and $F$ is strictly convex.  Let $T$ be an {\em into}
linear isometry from
$C_0(X,E)$ into $C_0(Y,F)$. Then there exist a continuous function  
$\varphi$
from a subset $Y_1$ of $Y$ onto $X$ and a continuous map $h$ from  
$Y_1$ into
$(B(E,F)$,$\SOT$)  such  that  for  all  $f$  in  $C_0(X,E)$,  
$$
Tf(y)  =  h(y)
(f(\varphi(y))),\quad \forall y \in Y_1.
$$ 
Moreover, $\|h(y)\|  = 1$, $\forall
y \in Y_1$, and for each $e$ in $E$ and $x$ in $X$, 
$$
\sup\{\|h(y)e\|: y \in
Y_1 \mbox{ and } \varphi(y) = x\} = \|e\|.
$$
Consequently,
$$
\|Tf\| = \|f\| = \|Tf_{|_{Y_1}}\|_\infty   \stackrel{\rm def}{=}
\sup_{y_1\in Y_1} \|Tf(y_1)\|.
$$
\end{thm}

It is easy to see that Jerison's  result \cite{Jer}
is a corollary  of Theorem \ref{thm}.
As indicated in \cite{Cam77},
there is a counter-example  in which the conclusion  of Theorem  
\ref{thm}
(in fact, even the one of Jerison \cite{Jer})  does not
hold while the assumption on strict convexity is not observed.
When $X$ and $Y$ are compact Hausdorff spaces, Theorem \ref{thm} reduces
to a result of Cambern \cite{Cam78}.  
It is plausible to think that Theorem \ref{thm} could be easily  
obtained from its
compact space version \cite{Cam78} by simply extending an into isometry
$T : C_0(X,E) \lrw C_0(Y,F)$ to an isometry from $C(X_\infty,E)$  
into  $C(Y_\infty,F)$
where $X_\infty$ (resp. $Y_\infty$) is the one-point  
compactification of $X$ (resp. $Y$).
However, an example in \cite[Example 9]{JW}
indicates that even in the simplest case $E=F={\mathbb R}$
there is an isometry from $C_0(X)$ into $C_0(Y)$ which cannot
be extended to an  isometry from $C(X_\infty)$ into $C(Y_\infty)$.  Thus
Theorem \ref{thm} cannot be obtained from the statement of the  
compact space version
directly.  
It is, however, possible to modify the argument in \cite{Cam78} to  
get a proof of
Theorem \ref{thm}.  The key is ``$F(x) = 0$ implies $AF(y)=0$'', in  
Cambern's
notation
\cite[Lemma 2]{Cam78}, which allows to define ``$A_y(e) = AF(y)$''  
where $F$ is any function
with $F(x) = e$.
Instead of going through the reasoning of Cambern once again, 
we present in the following an alternative approach based on the  
use of point evaluation
type functionals.
The technique of the proof we utilize here is
influenced  by those used in the scalar version  as appeared
in \cite{Hol} and \cite{Jar90}. 

    We would  like to take this opportunity  to express  our deep  
thanks  to
Cho-Ho Chu and Ka-Sing Lau for their encouragement.

\begin{proof}[Proof of  Theorem \ref{thm}]
    For a Banach  space $M$, we denote  by $U_M = \{m \in M:  
\|m\|\leq  1\}$
the closed  unit  ball, $S_M  =\{m\in  M: \| m\|=1\}$  the unit   
sphere, and
$M^\ast$  the Banach  dual space of $M$, respectively.  For $x$ in  
$X$, $y$
in $Y$, $\nu$ in $S_{E^\ast}$ and $\mu$ in $S_{F^\ast}$, we set
$$S_{x,\nu} = \{ f\in C_0(X,E): \nu(f(x))= \|f\|= 1\},$$
$$R_{y,\mu} = \{ g\in C_0(Y,F): \mu(g(y))= \|g\|= 1\},$$
$$Q_{x,\nu}  = \left\{\begin{array}{ll}  \{y\in  Y:  
T(S_{x,\nu})\subset  R_{y,
\mu}  \mbox{ for some } \mu \mbox{ in } S_{F^\ast}\},
& \mbox{if } S_{x,\nu}\not=
\emptyset,\\ \emptyset, &\mbox{if } S_{x,\mu}=\emptyset,\end{array}  
\right.$$
and $$Q_x=\bigcup_{\nu\in S_{E^\ast}} Q_{x,\nu}.$$

\bigskip\noindent {\sc Claim 1.} $Q_x \not= \emptyset$ for all $x$  
in $X$.

    Note that the product space $Y\times  U_{F^\ast}$  is a locally  
compact
Hausdorff  space.  Define  a linear  isometry  $\Psi$  from   
$C_0(Y,F)$ into
$C_0(Y\times  U_{F^\ast})$  by $$\Psi(g)(y,\mu)= \mu(g(y)).$$ Fix  
an $e$ in
$S_E$ and then a $\nu$ in $S_{E^\ast}$  such that $\nu(e)=  \|e\|=   
1$.  Then
$S_{x,\nu}\not= \emptyset$, $\forall x \in X$.  It now suffices to  
show that
$$\bigcap_{f\in  S_{x,\nu}} (\Psi(Tf))^{-1} \{1\} \not= \emptyset,\  
\ \forall
x\in X.$$ For each  $x$ in $X$, consider  $f_1,\ldots, f_n$  in  
$S_{x,\nu}$.
Let ${h=\sum_{i=1}^n f_i}$.  We have $\|h\|=n$ and thus there is a  
$y$ in $Y$
such that $\|Th(y)\|=n$.  So a $\mu$ in $S_{F^\ast}$  exists such that 
$n =
\mu(Th(y)) = \sum_{i=1}^n  \mu(Tf_i(y))$.  It then follows from  
$\|Tf_i(y)\| \leq
1$, $i=1,\ldots, n$, that  $\mu(Tf_i(y))$ $= 1$, $i=1, \ldots, n$,   
and  thus
$$(y,\mu)\in \bigcap_{i=1}^n  (\Psi(Tf_i)^{-1}(\{1\})  \not=  
\emptyset.$$ In
other  words, the family  $\{(\Psi(Tf))^{-1}  (\{1\}): f\in  
S_{x,\nu}\}$  of
compact sets has the finite intersection property.  Consequently,  
$Q_x \not=
\emptyset$.

\bigskip\noindent  {\sc Claim 2.} $Q_{x_1}\bigcap  Q_{x_2}  =  
\emptyset$  if
$x_1 \not= x_2$.

    Suppose  on the contrary  the  existence  of an $y$ in $Q_{x_1}  
 \bigcap
Q_{x_2}$.  Then there exist $\nu_1$ and $\nu_2$ in $S_{E^\ast}$   
and $\mu_1$
and $\mu_2$ in $S_{F^\ast}$ such that $$\mu_1(Tf(y))  =  
\nu_1(f(x_1)) = 1,\
\ \forall f \in S_{x_1,\nu_1}$$ and $$\mu_2(Tg(y))  = \nu_2(g(x_2))  
 = 1,\ \
\forall g \in S_{x_2,\nu_2}.$$ Let $U_1$ and $U_2$ be disjoint  
neighborhoods
of $x_1$ and $x_2$, respectively.  Choose $f_1$ in $S_{x_1,\nu_1}$   
and $f_2$
in $S_{x_2,\nu_2}$  such that $f_i$  is supported  in $U_i$, $i=1,  
2$.  Then
$\|f_1\pm  f_2\|  = 1$ implies  $\|T(f_1\pm  f_2)(y)\|  \leq 1$.   
In fact, the
inequalities  $2  =  2\|Tf_1(y)\|  =   
\|T(f_1+f_2)(y)+T(f_1-f_2)(y)\|   \leq
\|T(f_1+f_2)(y)\|  + \|T(f_1-f_2)(y)\|$  ensure that $\|T(f_1\pm  
f_2)(y)\| =
1$. By the strict convexity of $F$, we have $T(f_1+f_2)(y) =  
T(f_1-f_2)(y)$,
and thus a contraction that $Tf_2(y) = 0$!

    Let ${Y_1 = \bigcup_{x\in X} Q_x}$.  Define $\varphi : Y_1  
\rightarrow X$
such that $\varphi(y) = x$ if $y \in Q_x$.  For an $f$ in $C_0(X,E)$, we
denote $\coz f = \{ x \in X : f(x)\not= 0\}$ and $\supp f$ the  
closure of
$\coz f$ in $X$.  An argument similar to that in the proof of Claim  
2 will
give

    \bigskip\noindent{\sc Claim 3.} For each $f$ in $C_0(X,E)$,  
$\varphi(y)
\not\in \supp f$ implies $Tf(y) = 0$.

    \bigskip\noindent{\sc Claim 4.} $h(y)$ is well-defined and  
$\|h(y)\| =
1$ for all $y$ in $Y_1$.

    For each $y$ in $Y_1$, let $$J_y = \{f\in C_0(X,E) : \varphi(y)  
\not\in
\supp f \}$$ and $$K_y = \{f\in C_0(X,E) : f(\varphi(y)) = 0 \}.$$   
It is
not hard to see that $J_y$ is dense in $K_y$.  For $x$ in $X$  
(resp. $y$ in
$Y$), let $\delta_x$ (resp. $\delta_y$) be the point evaluation map
$\delta_x(f) = f(x)$ (resp. $\delta_y(g) =g(y)$) of $C_0(X,E)$ (resp.
$C_0(Y,F)$).  By Claim 3, $J_y \subset \ker(\delta_y\circ T)$ and thus
$\ker(\delta_{\varphi(y)}) = K_y \subset \ker(\delta_y\circ T)$.  Hence
there exists a linear operator $h(y)$ from $E$ into $F$ such that  
$$\delta_y
\circ T = h(y) \circ \delta_{\varphi(y)}.$$  In other words, for  
all $f$ in
$C_0(X,E)$, $$Tf(y) = h(y)(f(\varphi(y))).$$  For any $e$ in $E$,  
choose 
an $f$ in $C_0(X,E)$ such that $f(\varphi(y)) =
e$ and $\|f\| = \|e\|$.  Since $\|h(y)e\| = \|h(y)(f(\varphi(y)))\| =
\|Tf(y)\| \leq \|Tf\| = \|f\| =\|e\|$, we conclude that $\|h(y)\|  
\leq 1$.
In fact, it follows from the definition of $Y_1$ that $\|h(y)\| = 1$,
$\forall y \in Y_1$.  

The assertion that
for each $e$ in $E$ and $x$ in $X$, 
$\sup\{\|h(y)e\|: y \in
Y_1 \mbox{ and } \varphi(y) = x\} = \|e\|$
is obvious
if we pay attention to functions in the form of 
$f(w) = g(w)e$ where $g$ is a non-negative continuous
function on $X$ vanishing at infinity with maximum value $g(x) = 1$. 
Consequently, the norm identities $\|Tf\| = \|f\| =  
\|Tf_{|_{Y_1}}\|_\infty$
are established.

    \bigskip\noindent {\sc Claim 5.} $\varphi$ is continuous from  
$Y_1$ onto
$X$.

    Let $\{y_\lambda\}$ be a net convergent to $y$ in $Y_1$.  If $\{
\varphi(y_\lambda) \}$ does not converge to $\varphi(y)$, by  
passing to a
subnet if necessary, we assume it converges to an $x$ in $X_\infty = X
\bigcup \{ \infty \}$, the one-point compactification of $X$.  Let  
$U_1$ and
$U_2$ be disjoint neighborhoods of $x$ and $\varphi(y)$ in $X_\infty$,
respectively.  There exists a $\lambda_0$ such that  
$\varphi(y_\lambda) \in
U_1$, $\forall \lambda \geq \lambda_0$, and an $f$ in $C_0(X,E)$  
such that
$\coz f \subset U_2$ and $\|Tf(y)\| \not= 0$.  For  $\lambda \geq
\lambda_0$, $\varphi(y_\lambda) \not\in \supp f$.  By Claim 3,
$Tf(y_\lambda) = 0$, $\forall y_\lambda \geq \lambda_0$.  Thus
$\{Tf(y_\lambda)\}$ cannot converge to $Tf(y) \not= 0$, a  
contradiction.  Hence
$\varphi$ is continuous.

    \bigskip\noindent {\sc Claim 6.} $h : Y_1 \rightarrow  
(B(E,F),{\SOT})$ is continuous.

    Let $\{y_\lambda\}$  be a net convergent to $y$ in $Y_1$.  For  
$e$ in
$E$, $f$ in $C_0(X,E)$ exists such that $f(x) = e$ for all $x$ in a
neighborhood of $\varphi(y)$.  Since $\varphi$ is continuous, there is a
$\lambda_e$ such that for all $\lambda \geq \lambda_e$,  
$\|h(y_\lambda)e -
h(y)e\| = \|h(y_\lambda)f(\varphi(y_\lambda)) - h(y)f(\varphi(y))\| =
\|Tf(y_\lambda) - Tf(y)\|.$  Since $\{Tf(y_\lambda)\}$ converges to
$Tf(y)$, the claim is thus verified.  The proof is complete. 
\end{proof}

To end this note, we would like to remark that $Y_1$ can be neither  
open nor closed,
and $h(y)$ need not be an isometry in  general for $y$ in $Y_1$ as  
pointed out by
an example in \cite{Cam78}.  

\bibliographystyle{plain}

    \end{document}